\documentclass[10pt,openany,leqno]{smfart}

\usepackage[english, francais]{babel}
\usepackage{smfthm}
\usepackage{amsmath,amsthm,amsfonts,amssymb,amscd,url}
\usepackage{graphics}

\input{xy} \xyoption{all}
\usepackage[latin1]{inputenc}

\usepackage{enumerate}
\usepackage{hyperref}
\usepackage{epsfig} 
\input{xy} \xyoption{all}

\begin{document}
\def\Diff{\text{Diff}}
\def\Max{\text{max}}
\def\R{\mathbb R}
\def\N{\mathbb N}
\def\Z{\mathbb Z}
\def\C{\mathbb C}
\def\T{\mathbb T}
\def\a{{\underline a}}
\def\b{{\underline b}}
\def\c{{\underline c}}
\def\Log{\text{log}}
\def\loc{\text{loc}}
\def\inta{\text{int }}
\def\det{\text{det}}
\def\exp{\text{exp}}
\def\Re{\text{Re}}
\def\lip{\text{Lip}}
\def\leb{\text{Leb}}
\def\dom{\text{Dom}}
\def\diam{\text{diam}\:}
\def\supp{\text{supp}\:}
\newcommand{\ovfork}{{\overline{\pitchfork}}}
\newcommand{\ovforki}{{\overline{\pitchfork}_{I}}}
\newcommand{\Tfork}{{\cap\!\!\!\!^\mathrm{T}}}
\newcommand{\whforki}{{\widehat{\pitchfork}_{I}}}
\newcommand{\marginal}[1]{\marginpar{{\scriptsize {#1}}}}
\def\sR{{\mathfrak R}}
\def\sM{{\mathfrak M}}
\def\sA{{\mathfrak A}}
\def\sB{{\mathfrak B}}
\def\sY{{\mathfrak Y}}
\def\sE{{\mathfrak E}}
\def\sP{{\mathfrak P}}
\def\sG{{\mathfrak G}}
\def\sa{{\mathfrak a}}
\def\sb{{\mathfrak b}}
\def\sc{{\mathfrak c}}
\def\se{{\mathfrak e}}
\def\sg{{\mathfrak g}}
\def\sd{{\mathfrak d}}
\def\sr{{\mathfrak {r}}}
\def\ss{{\mathfrak {s}}}
\def\st{{\mathfrak {t}}}
\def\spp{{\mathfrak {p}}}
\def\avv{ \overrightarrow }
\def\arr{ \overleftarrow }
\theoremstyle{plain}
\newenvironment{prob}{\begin{enonce}{Problem}}{\end{enonce}}




\title{Strong regularity}

\author{Pierre Berger and Jean-Christophe Yoccoz}

\date{Submitted 06/16/2014, accepted 01/21/2018}
\begin{otherlanguage}{english}
\maketitle

\section{Uniformly hyperbolic dynamical systems}

The theory of \emph{uniformly hyperbolic dynamical systems} was constructed in the 
1960's under the dual leadership of Smale in the USA and Anosov and Sinai in the Soviet Union.
It is nowadays  almost complete. It encompasses various examples \cite{Sm}: expanding maps, horseshoes, solenoid maps, Plykin attractors,  Anosov maps and DA, all of which are \emph{basic pieces}. 

\par

We recall standard definitions. Let $f$ be a $C^1$-diffeomorphism $f$ of a finite dimensional manifold $M$. A compact $f$-invariant subset $\Lambda \subset M$ is 
\emph{uniformly hyperbolic} if the restriction to $\Lambda$ of the tangent bundle $TM$
splits into two continuous invariant subbundles
\[TM|\Lambda = E^s\oplus E^u,\]
$E^s$ being  uniformly contracted and $E^u$ being uniformly expanded. 

Then for every $z\in \Lambda$, the sets 
\[W^s(z)= \{z'\in M:\; \lim_{ n\to+\infty} d(f^n(z),f^n(z'))= 0\},\]
\[W^u(z)= \{z'\in M:\; \lim_{ n\to-\infty} d(f^n(z),f^n(z'))= 0\}\]
are called \emph{the stable and unstable manifolds} of $z$. They are immersed manifolds tangent at $z$ to respectively $E^s(z)$ and $E^u(z)$.

The  \emph{$\epsilon$-local  stable manifold} $ W^s_\epsilon(z)$ of $z$ is the connected component of $z$ in the intersection of $W^s(z)$ with a $\epsilon$-neighborhood  of $z$. The \emph{$\epsilon$-local  unstable manifold} $ W^u_\epsilon(z)$ is defined likewise.
%



\begin{defi} 
A \emph{basic set}  is a compact, $f$-invariant, uniformly hyperbolic  set $\Lambda$ which is transitive and \emph{locally maximal}: there exists   a neighborhood $N$ of $\Lambda$ such that $\Lambda = \cap_{n\in \Z} f^n(N)$. A basic set is an \emph{attractor} if the neighborhood $N$ can be chosen in such a way that $\Lambda = \cap_{n \geq 0} f^n(N)$. Such a basic set contains the unstable manifolds of its points.
\end{defi}

 A diffeomorphism whose nonwandering set is a finite union of disjoint  basic sets is called  \emph{uniformly hyperbolic} or \emph{Axiom A}.
 
 \bigskip

Such diffeomorphisms enjoy nice properties, which are proved in \cite{Sm} and the references therein.
  
\subsubsection*{SRB and physical measure}

Let $\alpha >0$, and let $\Lambda$ be an attracting basic set for a \linebreak $C^{1+\alpha}$-diffeomorphism $f$. Then there exists a unique invariant, ergodic probability $\mu$ supported on $\Lambda$ such that its conditional measures, with respect to any measurable partition of $\Lambda$ into plaques of unstable manifolds, are absolutely continuous with respect to the Lebesgue measure class (on unstable manifolds). Such a probability is called \emph{SRB} (for Sinai-Ruelle-Bowen). It turns out that a SRB -measure is \emph{physical}: the Lebesgue measure of its basin $B(\mu)$
\begin{equation}\tag{$\mathcal B$} B(\mu) = \{z\in M: \; \frac1n \sum_{i < n} \delta_{f^i(x)} \rightharpoonup \mu\},\end{equation} 
 is positive. Actually, up to a set of Lebesgue measure $0$, $B(\mu)$ is equal to the topological basin of $\Lambda$, i.e the set of points attracted by $\Lambda$.

\subsubsection*{Persistence} A basic set $\Lambda$ for a $C^1$-diffeomorphism $f$ is \emph {persistent}: every  $C^1$-perturbation $f'$ of $f$ leaves invariant a basic set $\Lambda'$ which is homeomorphic to $\Lambda$, via a homeomorphism which conjugates the dynamics $f|\Lambda$ and $f'|\Lambda'$.


\subsubsection*{Coding} A basic set $\Lambda$ for a $C^1$-diffeomorphism $f$ admits a (finite) Markov partition. This implies that its dynamics is semi-conjugated with a subshift of finite type. The semi-conjugacy is 1-1 on a generic set. Its lack of injectivity is itself coded by subshifts of finite type of smaller topological entropy. This enables to study efficiently all the invariant measures of $\Lambda$, the distribution of its periodic points, the existence and uniqueness of the maximal entropy measure, and if $f$ is $C^{1+\alpha}$, the Gibbs measures which are related to the geometry of $\Lambda$.

\subsection{End of Smale's program}
Smale wished to prove the density of Axiom A in the space of $C^r$-diffeomorphisms. 
In higher dimensions,  obstructions were soon discovered by Shub \cite{Sh71}. For surfaces Newhouse showed the non-density of Axiom A diffeomorphisms  for $r\ge 2$: he constructed  robust tangencies between stable and unstable manifolds of a thick horseshoe \cite{Newhouse}. Numerical studies by Lorenz \cite{Lorenz} and  H\'enon \cite{Henon} explored dynamical systems with hyperbolic features that did not fit in the uniformly hyperbolic theory.
In order to include many examples such as the H\'enon one, the \emph{non-uniform hyperbolic theory} is still under construction.

\section{Non-uniformly hyperbolic dynamical systems}
\subsection{ Pesin theory}
The natural setting for non-uniform hyperbolicity is Pesin theory \cite{BP06,LY}, from which we recall some basic concepts. We first consider the simpler settings of invertible dynamics.

\par

Let $f$ be a $C^{1+\alpha}$-diffeomorphism  (for some $\alpha >0$) of a compact manifold $M$ and let $\mu$ be an ergodic $f$-invariant probability measure on $M$.
The Oseledets multiplicative ergodic theorem produces  Lyapunov exponents (w.r.t. $\mu$)
for the tangent cocycle of $f$, and an associated $\mu$-a.e $f$-invariant splitting of the tangent bundle into characteristic subbundles.

Denote by $E^s(z)$ (resp. $E^u(z)$) the sum of the the characteristic subspaces associated to the negative (resp. positive) Lyapunov exponents.

The \emph{  stable and unstable Pesin manifolds} are defined respectively for $\mu$-a.e. $z$ by

\[W^s(z)= \{z'\in M:\limsup_{n\to+\infty} \frac1n \log d(f^n(z),f^n(z'))<0\},\]
\[ W^u(z)= \{z'\in M: \liminf_{n\to-\infty} \frac1n \log d(f^n(z),f^n(z'))>0\}.\]

They are immersed manifolds through $z$ tangent respectively at $z$ to $E^s(z)$ and $E^u(z)$.

The measure $\mu$ is \emph{hyperbolic} if $0$ is not a Lyapunov exponent w.r.t. 
$\mu$.
Every invariant ergodic measure, which is supported on a uniformly hyperbolic compact invariant set, is hyperbolic. 

\subsubsection*{SRB, physical measures}
An invariant ergodic measure $\mu$ is  \emph {SRB} if the largest Lyapunov exponent is positive  and the conditional measures  of $\mu$ w.r.t. a measurable partition into plaques of  unstable manifolds are $\mu$-a.s. absolutely continuous w.r.t. the Lebesgue class (on unstable manifolds). 
When $\mu$ is SRB and hyperbolic, it is also \emph{ physical}: its basin has positive Lebesgue measure.

\par
The paper \cite{Y98} provides a general setting where appropriate hyperbolicity hypotheses allow to construct hyperbolic SRB measures with nice statistical properties.

\subsubsection*{Coding}
Let $\mu$ be a $f$-invariant ergodic hyperbolic SRB measure. Then there is a partition mod.$0$ of $M$ into finitely many disjoint subsets $\Lambda_1,\ldots, \Lambda_k$, which are cyclically permuted by $f$ and such that the restriction $f^k_{| \Lambda_1}$ is metrically conjugated to a Bernoulli automorphism.

Of a rather different flavor is Sarig's recent work \cite{Sa13}. For a $C^{1+\alpha}$-diffeomorphism of a compact surface
of positive topological entropy and any $\chi >0$, he constructs a countable Markov partition for an invariant set which has full measure w.r.t. any ergodic invariant measure with metric entropy $>\chi$. The semi-conjugacy associated to this Markov partition is finite-to-one.



\subsubsection*{Non-invertible dynamics}

One should distinguish between the non-uniformly expanding case and the case of general endomorphisms.

\smallskip

In the first setting, a SRB measure is simply an ergodic invariant measure whose Lyapunov exponents are all  positive and which is absolutely continuous.

\smallskip

Defining appropriately unstable manifolds and SRB measures for general endomorphisms is more delicate. One has typically to introduce the inverse limit where the endomorphism becomes invertible.

\subsection{Case studies}

The paradigmatic examples in low dimension can be summarized by the following table:  

\begin{center}
	\begin{tabular}{|c|c|}
	\hline
		Uniformly hyperbolic & Non-uniformly hyperbolic\\
		\hline
		Expanding maps of the circle      &  Jakobson's Theorem \\
		Conformal expanding maps of complex tori & Rees' Theorem \\
		Attractors (Solenoid, DA, Plykin...) &  Benedicks-Carleson's Theorem \\
		Horseshoes & Non-uniformly hyperbolic horseshoes \\
		Anosov diffeomorphisms & Standard map ?\\
		\hline
	\end{tabular}
\end{center}

Let us recall what these theorems state, and the correspondence given by the lines of the table.


Expanding maps of the circle may be considered as the simplest case of uniformly hyperbolic dynamics. The Chebychev quadratic polynomial $P_{-2}(x) := x^2 -2$
on the invariant interval $[-2,2]$ has a critical point at $0$, but it is still semi-conjugated to the doubling map $\theta \mapsto 2\theta$ on the circle (through $x = 2 \cos 2 \pi \theta$). For $a \in [-2, -1]$, the quadratic polynomial $P_a(x):=x^2 +a$ leaves invariant the interval $[P_a(0), P^2_a(0)]$ which contains the critical point $0$.

\begin{theo}[Jakobson \cite{Ja81}] There exists a   set $\Lambda\subset [-2,-1]$ of positive Lebesgue measure  such that for every $a\in \Lambda$ the map $P(x)=x^2+a$ leaves invariant an ergodic, hyperbolic measure which is equivalent to the Lebesgue measure on $[P_a(0), P^2_a(0)]$.

\end{theo}
Actually the set $\Lambda$ is nowhere dense. Indeed the set of $a \in \R$ such that $P_a$ is  axiom A is open and dense \cite{GS97, Ly97}.

\bigskip

Let $L$ be a lattice in $\C$ and let $c$ be a complex number such that $|c| >1$ and $c L \subset L$. Then the homothety $ z \mapsto cz$ induces an expanding map of the complex torus $\C /L$. The Weierstrass function associated to the lattice $L$ defines 
a ramified covering of degree $2$ from $\C /L$ onto the Riemann sphere which is a semi-conjugacy  from this expanding map to a rational map of degree $|c|^2$ called a \emph{Lattes map}. For any $d \ge 2$, the set ${\rm Rat}_d$ of rational maps of degree $d$ is naturally parametrized by an open subset of $\mathbb P(\C^{2d+2})$.

\begin{theo}[Rees \cite{Re86}] For every $d\ge 2$, there exists a subset $\Lambda \subset {\rm Rat}_d$ of positive Lebesgue measure such that every  map $R\in \Lambda$ leaves invariant an ergodic hyperbolic probability measure which is equivalent to the Lebesgue measure on the Riemann sphere.
\end{theo}

For rational maps in $\Lambda$, the Julia set is equal to the Riemann sphere.
On the other hand, a conjecture of Fatou \cite{Mi06} claims that the set of rational maps which satisfy Axiom A is open and dense in ${\rm Rat}_d$. The  restriction of such maps to their  Julia set is uniformly expanding. 
For such maps, the Hausdorff dimension of the Julia set is  smaller than $2$. 

\bigskip

The (real) H\'enon family is the $2$-parameter family of polynomial diffeomorphisms of the plane defined for $a,b \in \R$, $b \ne 0$ by 

\[h_{a\,b}  (x,y) = (x^2+a+y,-bx)\]

Observe that $h_{a\,b}$ has constant Jacobian equal to $b$.
For small $|b| $, there exists an interval $J(b)$ close to $[-2, -1]$ such that, for $a \in J(b)$, the H\'enon map $h_{a\,b}$ has the following properties
\begin{itemize}
\item $h_{a\,b}$  has two fixed points; both are hyperbolic saddle points, one, called $\beta$  with positive unstable eigenvalue, the other , called $\alpha$, with negative unstable eigenvalue;
\item there is a trapping open region $B$ satisfying $h_{a\,b}(B) \Subset B$ which contains $\alpha$  (and therefore also its unstable manifold ).
\end{itemize}

H\'enon \cite{Henon} investigated numerically the behavior of orbits starting in $B$ for $b=-0.3$, $a=-1.4$. Such orbits apparently converged to a \textquotedblleft strange attractor \textquotedblright. 

\begin{theo}[Benedicks-Carleson \cite{BC2}] For every $b<0$ close enough to $0$, there exists a set $\Lambda_b\subset J(b)$ of positive Lebesgue measure, such that for every $a \in \Lambda_b$, the maximal invariant set $\bigcap_{n\ge 0} h^n_{a\,b}(B)$ is equal to the closure of the unstable manifold $W^u(\alpha)$ and contains a dense orbit along which the derivatives of iterates grow exponentially fast.
\end{theo}

An easy topological argument ensures that this maximal invariant set is never uniformly hyperbolic.
Later Benedicks-Young \cite{BY} showed that for every such parameters $a\in \Lambda_b$ the Hénon map $h_{a\,b}$ leaves invariant an ergodic hyperbolic SRB  measure. Such a measure is physical. Benedicks-Viana \cite{BV01} actually proved that the basin of this measure has full Lebesgue measure in the trapping region $B$. 

From \cite{Ur95}, every  $a \in \Lambda_b$ is accumulated by parameter intervals exhibiting Newhouse phenomenon:  for generic parameters in these intervals, $h_{a\,b}$ has infinitely many periodic sinks in $B$. In particular, the set $\Lambda_b$ is nowhere dense.



\bigskip

The starting point in \cite{PY09} is a smooth diffeomorphism of a surface $M$ having a horseshoe \footnote{A horseshoe is an infinite basic set of saddle type.} $K$. It is assumed that there exist distinct fixed points $p_s,p_u \in K$ and $q \in M$ such that 
$W^s(p_s)$ and $W^u(p_u)$ have at $q$ a quadratic heteroclinic tangency which is an isolated point of $W^s(K) \cap W^u(K)$. The authors consider a one-parameter family $(f_t)$ unfolding the tangency and study the maximal $f_t$-invariant set $L_t$
in a neighborhood of the union of $K$ with the orbit of $q$. Writing $d_s,d_u$ for the transverse Hausdorff dimensions of $W^s(K) ,\; W^u(K)$ respectively, it was shown previously \cite{PT93} that $L_t$ is a horseshoe for most $t$ when $d_s + d_u <1$. By \cite{MY10} this is no longer true when $d_s + d_u >1$. However, when $d_s + d_u $ is only slightly larger\footnote{The exact condition is  $(d_s+d_u)^2+(\max\{d_s,d_u\})^2<d_s+d_u+\max\{d_s,d_u\}$ .} than $1$, some dynamical and geometric information on $L_t$ is obtained in \cite{PY09} for most values of $t$: in particular, both the stable and unstable sets for $L_t$ have Lebesgue measure $0$, and an ergodic hyperbolic $f_t$-invariant probability measure supported on $L_t$ with geometric content is constructed.

\bigskip

The two papers in this volume are related to these case studies.

\medskip

In \cite{Y95}, a proof of Jakobson's theorem is given. The main  ingredient is the concept of \emph {strong regularity} (explained below).

\medskip

In  \cite{B11}, a class of endomorphisms of the plane containing the H\'enon family is considered. Given any map $B\in C^2(\R^3,\R^2)$ with small $C^2$-uniform norm, one studies the one-parameter family
\[f_{a,B}(x,y)= (x^2+a,0)+B(x,y,a).\]
It is shown that there exists a set $\Lambda_B \subset \R$ of positive Lebesgue measure such that, for any $a \in \Lambda_B$, $f_{a,B}$ has an invariant ergodic hyperbolic physical SRB measure. The proof is based on an appropriate generalization of strong regularity.





%

\subsection{	Open problems}

Linear Anosov diffeomorphisms of $\T^2$ are area-preserving and uniformly hyperbolic. In the conservative setting, a very natural case study to consider is the Chirikov-Taylor standard map family. This is a one-parameter family of area-preserving diffeomorphisms of $\T^2$ defined for $a \in \R$ by
\[ S_a(x,y) = (2x -y + a \sin 2 \pi x, x).\]
One form of a conjecture of Sinai ( \cite{Si94} P.144) about this family is 
\begin{conj}
There exists a set $\Lambda\subset\R$ of positive Lebesgue measure such that, for $a \in \Lambda$,  the Lebesgue measure on $\T^2$ is ergodic and hyperbolic for $S_a$. 
\end{conj}

For such parameters, the map $S_a$ cannot have any of the invariant curves produced by KAM-theory. In particular, $a$ cannot be too small.

This conjecture is still completely open despite intense efforts. A weak argument in favor of this conjecture is that, when $a$ is large, the maximal invariant set in the complement of an appropriate neighborhood of the critical lines $\{x= \pm 1/4 \}$ is a uniformly hyperbolic horseshoe of dimension close to 2 \cite{Du94, BeCa13}.


Actually, a large  Hausdorff dimension of the invariant sets under consideration appears to be a major difficulty on the way to prove non-uniform hyperbolicity. 

For the parameters considered in  \cite{BC2} and subsequent papers, the Hausdorff dimension of the H\'enon attractor is \emph {a priori}  close to $1$. On the other hand, numerical studies \cite{RHO80} of the values $a=-1.4$, $b=-0.3$ considered by H\'enon
indicate an (eventual) attractor of Hausdorff dimension $1.261 \pm 0.003$.

\begin{prob}\label{P1}
For every $d<2$, find an open set of smooth families  $(f_t)_t$ of smooth diffeomorphisms  of $\R^2$ such that, with positive probability on the parameter,  $f_t$ leaves invariant an ergodic hyperbolic SRB probability measure whose support has dimension at least $d$.   
\end{prob} 

One should also recall that  Carleson conjectured \cite[p. 1246 l.18]{Ca91} that proving non-uniform hyperbolicity
[or only the weaker conclusion of  \cite{BC2}] ``for a particular parameter value is in some rigorous sense undecidable''.

A similar problem, in the setting of non-uniformly hyperbolic horseshoes, is 

\begin{prob}\label{P2}
Prove the conclusions of \cite{PY09} for an initial horseshoe $K$ of transverse Hausdorff dimensions $d_s,d_u$ satisfying 
\[d_s+d_u>3/2.\]    
\end{prob}

Even the non-uniformly expanding case is still incomplete, since it considers only the case of real or complex dimension 1. A positive answer to the following problem would be a 2-dimensional generalization of Jakobson's Theorem for perturbation of the product dynamics:
\[P_a\times P_a\colon (x,y)\mapsto (x^2+a,y^2+a).\]

\begin{prob}
Does there exist an open set of $1$-parameter smooth families $(f_a)$ of endomorphisms of the plane, accumulating on $(P_a \times P_a)_a$, with the following property: with positive probability on the parameter, $f_a$ leaves invariant
an ergodic absolutely continuous invariant measure with two positive Lyapunov exponents.
\end{prob}


%
%
%

\section{Proving non uniform hyperbolicity}

There are now many proofs of both Jakobson's theorem and Benedicks-Carleson's theorem. Broadly speaking, they rely either on a binding approach, pioneered by Benedicks-Carleson, or on a strong regularity approach, closer to Jakobson's original proof \cite{Ja81} and to \cite{R88}.
Both papers in this volume follow the second approach. 

In both approaches, the study of the $2$-dimensional setting depends very much on the $1$-dimensional case.

We now explain some of the differences between the two methods.

\subsection{ The binding approach  for quadratic maps}

Benedicks-Carleson proved Jakobson's theorem by focusing on the expansion of the post-critical orbit. There are many proofs in this spirit \cite{CE80, BC1, Ts932, Ts93, Lu00}.

One actually proves the existence of a  set $\Lambda \subset \R$ of positive Lebesgue measure such that, for $a \in \Lambda$, the quadratic map $P_a(x)=x^2+a$ satisfies the  Collet-Eckmann condition: 
\[\liminf_{+\infty} \frac1n \log \|DP^n(a)\|>0.\] 
This property implies the existence of an absolutely continuous ergodic invariant measure with positive Lyapunov exponent\cite{CE83}.

One starts with a parameter $a_0$   such that the critical value $a_0$ of $P_{a_0}$ belongs to a repulsive periodic cycle. Then, there exists $\lambda>1$ so that 
\begin{itemize}
\item[$(i)$] $DP_{a_0}^n(a_0)>\lambda^n$ for every  large $n$,
\item[$(ii)$] for every $\delta>0$, the map $P_{a_0}$ is $\lambda$-expanding on the complement of $[-\delta,\delta]$ (for an adapted metric).  
\end{itemize}

Then for every  large $M$,  for every $a$ close to $a_0$ the post-critical orbit $(P^n_a(a))_{n\le M}$ is close to $(P_{a_0}^n(a_0))_{n\le M}$ and so has a similar expansion. At the next iterations $N=M+1$, there are three possibilities:
\begin{itemize}
\item[$(a)$] either $P_{a}^{N}(a)$ is not in $(-\delta,\delta)$ and so the expansion will continue by $(i)$, 
\item[$(b)$] or $P_{a}^{N}(a)$ is in $(-\delta,\delta)$ but is not too close to $0$; then there exists an integer $k<N$, called \emph{the binding time}, such that the orbits $P_{a}^{N+i}(a)$ and $P_{a}^{i}(0)$ remain close for $i\le k$ and separate for $i=k+1$. The expansion of $(DP_{a}^{i}(a))_{i <k}$ is transferred to $(DP_{a}^{i}(P^{N+1}_a(a)))_{i  < k}$ . The  logarithmic contraction at time $N$, equal to $\log |DP_a (P_a^N (a))|$, is only roughly half the logarithmic expansion during the binding period $\log |DP^{k-1}_a (P_a^{N+1} (a))|$.
\item[$(c)$] or $P_{a}^{N}(a)$ is so close to $0$ that $(b)$ does not hold.
\end{itemize}
Cases $(a)$ and $(b)$ are allowed. Case $(c)$ is excluded in the parameter selection by removing the parameter $a$ for which this occurs.
Then we can redo the same alternative with $N\rightarrow N+1$ in case $(a)$ and $N\rightarrow N+k$ in case $(b)$.

In case $(b)$, roughly half of the original transferred logarithmic expansion is lost in the binding process. Therefore the Collet-Eckmann condition will not be satisfied if too much time is spent in iterated binding periods. To avoid this, it is asked that:

\begin{itemize}
\item[$(H_N)$] the total length of all the  binding periods before $N$ is small with respect to $N$.
\end{itemize}
Actually, when appropriately formulated, the condition $(H_N)$ implies that case $(c)$ above does not hold. Hence if $(H_N)$ holds for every $N$, the map is Collet-Eckmann. 

To perform the parameter selection, we look at maximal \emph{critical curves} $\gamma=(P_a^N(a))_{a\in \mathcal I}$ so that:
\begin{itemize} 
\item[($P_1$)] Condition $(H_n)$ holds for  every $a\in \mathcal I$ and for every $n\le N$;
\item[($P_2$)] the binding periods in $[0,N]$ are the same for every $a\in \mathcal I$, and the integer $N$ is not part of a binding period;
\item[($P_3$)] the length of the curve $\gamma$ is bounded from below by some uniform constant.
\end{itemize}

Such a curve is split into different pieces according to which scenario holds at time $N+1$. Pieces corresponding to scenario $(a)$ are iterated once. Pieces corresponding to scenario $(c)$ (or to scenario $(b)$, with a binding time $k$ too long to satisfy $(H_{N+k})$) are discarded. The other pieces are iterated untill the end of the corresponding binding period. These new critical curves satisfy $(P_1)$ and $(P_2)$. Property $(P3)$ is also satisfied, except for some boundary effects that are easily taken care of.

\medskip


A large deviation argument, relying on property $(P3)$, shows that the Lebesgue measure of the remaining parameters is positive (actually, a large proportion of the length of the starting parameter interval).


\subsection{The binding approach for H\'enon family}

There are many proofs in this spirit \cite{BC2, MV93, WY01,VL03, YW08, Ta11}. 

\medskip

A major difficulty of the $2$-dimensional setting is that critical points are not defined beforehand, and will only be well-defined for good parameters.

Call a curve \emph{flat} if it is $C^2$-close to a segment  of $\R\times \{0\}$. Roughly speaking, given a flat segment $\gamma \subset W^u(\alpha)$ going across the critical strip $\{ \vert x \vert \leq \delta \}$, a critical point on $\gamma$ should be  a point of $\gamma$
 such that the vertical tangent vector is exponentially dilated under positive iteration, while the tangent vector to $\gamma$ is exponentially contracted.
 
 \smallskip
 
 In the inductive construction of good parameters, only $N$ iterations of the H\'enon map are considered at a given stage. Under the appropriate induction hypotheses, one defines an approximate critical set $\mathcal C_N$. This is a finite set whose cardinality is exponentially large with $N$. Each point of $\mathcal C_N$ lies on a flat segment contained in $h_{a\,b}^{\theta N} (W^u_{loc}(\alpha))$, with $\theta\sim \vert \log \vert b \vert \vert^{-1}$. 
 
 The main problem of the induction step is to extend the exponential dilation along the finitely many critical orbits beyond time $N$. As in the $1$-dimensional case, this is automatic when the critical orbit at time $N$ lies outside of the critical strip. On the other hand, when the critical orbit at time $N$ returns to a point $z_N$ of the critical strip, one has to find, after excluding inadequate parameters, a \emph{binding} critical point $\tilde z_0$ whose initial expansion will be transferred (at some cost) to the orbit of $z_N$.
 It is here important that $z_N$ should be in \emph{tangential position}, i.e much closer to the flat segment containing $\tilde z_0$ than to $\tilde z_0$ itself.
 
 \smallskip
 
 To prove that the set of non-excluded parameters (at the end of the induction process) has positive Lebesgue measure, one has to investigate carefully how the whole structure of approximate critical points, analytical estimates and binding relationships survives through parameter deformation. This is certainly the  trickiest part of the method.

\subsection{Puzzles and parapuzzles}
Puzzles and parapuzzles are combinatorial structures which were first introduced in $1$-dimensional complex dynamics to study the local connectivity of Julia sets and the Mandelbrot set  \cite{Hu91,Mi92}. In real $1$-dimensional dynamics, they were instrumental in the proof that almost every quadratic map satisfies  either axiom A or the Collet-Eckmann condition \cite{L02, AM03}.

\par
For real Julia sets of real quadratic maps, puzzle pieces are defined as follows. Let $a$ be a parameter in $[-2,-1]$. Then the quadratic polynomial $P_a$ has two fixed points $\alpha,\beta$, both repelling, denoted so that $-\beta < \alpha < -\alpha < \beta$. The real Julia set is equal to $[-\beta, \beta]$. For $n\geq 0$, the \emph{puzzle pieces of order $n$} are the closures of the connected components of $[-\beta, \beta] \setminus P_a^{-n}(\{\alpha,-\alpha \})$.

\par
Puzzle pieces of successive orders are related in two fundamental ways: a puzzle piece of order $n$ is contained in a puzzle piece of order $n-1$, and its image is contained in a puzzle piece of order $n-1$. The combinatorics of the partition by puzzle pieces of a given order depend on the sequence of nested puzzle pieces containing the critical value. This leads to a sequence of partitions of parameter space into \emph{parapuzzle} pieces. It is a general rule of thumb that, assuming a mild level of hyperbolicity, the combinatorics and geometry of parapuzzle pieces around a given parameter $a$ are closely related to the combinatorics and geometry of puzzle pieces for $P_a$ around the critical value. 

\subsection{The strong regularity approach  for quadratic maps}\label{SR1}

 Let $a$ be a parameter in $[-2,-1]$. A \emph{regular interval} is a puzzle piece of some order $n>0$ which is sent diffeomorphically onto $A := [\alpha, -\alpha]$ by $P_a^n$. One also asks that the corresponding inverse branch extends to a fixed neighborhood of $A$, which insures a control of the distortion. The parameter $a$ is  \emph{regular} if the measure of the set of points in $A$ which are not contained in a regular interval of order $\leq n$ is exponentially small with $n$. A classical argument shows that regular parameters satisfy the conclusions of Jakobson's theorem.

\par
To prove that the set of regular parameters has positive Lebesgue measure, one considers a more restrictive condition called \emph{strong regularity}. Assume that the parameter is close to the Chebychev value $a_0:= -2$. Then the return time $M$ of the critical point to $A$ is large. Moreover, the complement in $A$ of a neighborhood of $0$ of approximate size $2^{-M}$ is covered by finitely many regular intervals of order $<M$, which are called \emph{simple}. The parameter $a$ is called \emph{strongly regular} if
\begin{itemize}
\item[$(\star)$] there exists a sequence of regular intervals $(I_j)_{j>0}$ of order $(n_j)_j$ such that   $P_a^{M+n_1+\cdots+ n_{j-1}}(a)\in I_{j}$ for all $j>0$;
\item[$(\diamondsuit )$] most $I_j$ are simple in the sense that $\sum_{i\le j: I_i\text{ is not simple} }n_i<\!\!< \sum_{i\le j} n_i$ for all $j>0$.\label{SRcondition}
\end{itemize}

The most delicate part of the proof is to establish, through a careful analysis of the puzzle structures, that strongly regular parameters are regular. Then one is able to transfer the exponential regularity estimate from puzzles in phase space to parapuzzles in parameter space. Finally, one concludes through a large deviation argument that the set of strongly regular parameters has positive Lebesgue measure.

\subsection{ The strong regularity approach for H\'enon family}\label{SR2}

We assume that the $C^2$-norm $b$ of $(B_a)_a$ is small, and we put $\theta= 1/|\log b|$. The hyperbolic fixed point $(\alpha,0)$ of $h_{a\, 0}$ persists as a fixed point $A$ for $f=h_{a\, B}: (x,y)\mapsto (x^2+a, 0)+B_a(x,y)$. 
Up to a conjugacy we can assume that the vertical boundary $\partial^s Y_\se:= \{\alpha, -\alpha\}\times [-\theta, \theta]$ of $ Y_\se:=[\alpha, -\alpha]\times [-\theta, \theta]$ is such that $\{\alpha\}\times [-\theta, \theta]$ is a local stable manifold  of $f$ and that $\{-\alpha\}\times [-\theta, \theta]$ is sent by $f$ into $\{\alpha\}\times [-\theta, \theta]$. 

A \emph{box} $Y$ is a subset of $\R\times [-\theta, \theta]$ bounded by two arcs $\partial^s Y$ of the stable manifold of $A$ which are $C^2$-close to be vertical and with endpoints in  $\R\times \{-\theta, \theta\}$. A \emph{piece} $(Y,n)$ is  the data of a box $Y$ and an integer $n$ such that  $f^n(Y)\subset Y_\se$ and any nearly horizontal vector $u$ pointed at $z\in Y$ verifies conditions on uniform expansion on its $n$-first iterates.
A \emph{puzzle piece} $(Y,n)$  is a piece such that $f^n(\partial^s Y)\subset \partial^s Y_\se$. \medskip

The simple puzzle pieces of the one-dimensional Chebichev map survive for the studied 2-dimensional endomorphisms to form a set of puzzle pieces in $Y_\se$  denoted by $\{(Y_\ss, n_\ss): \ss\in \sA_0\}$, where $\sA_0$ is a set of $2M-2$ symbols and $M$ is the first return time of $(a,0)$ in  $Y_\se$. The complement of the union of these simple pieces in $Y_\se$ is a box $Y_\boxdot$ of approximate width $2^{-M}$.  
\medskip
 
Similarly to \cite{PY09}, we define two operations on the pieces.  The first is the $\star$-\emph{product} between two pieces $(Y,n)$ and $(Y', n')$ in $Y_\se$, equal to $(Y,n)\star (Y',n')=(f^{-n}(Y')\cap Y, n+n')$. Whenever the constructed set is not included in an arc of $W^s(A)$, the pair obtained is a piece and the product is called \emph{admissible}. 

The second operations are the \emph{parabolic products} denoted by $\boxdot_\pm$ for $\pm\in \{-,+\}$. They enable to construct pieces with iterations visiting the region $Y_\boxdot$. They are considered when a puzzle piece $(Y,n)$ satisfies that $f^n(Y) $ is folded by $f^{M+1}|Y_\boxdot$ at a puzzle piece $Y'$ but not too close to $\partial^s Y'$. If $(Y'',n'')$ is a piece such that $Y''\supset Y'$, the left component of $\partial^s Y'$ and $\partial^s Y''$ are different, and the orders $n, n',n''$ satisfy some linear bounds, then this triplet of pieces is called \emph{pre-admissible} for the parabolic product. The preimage of $(f^{M+1}|Y_\boxdot)\circ (f^n|Y)$ of $Y''\setminus int\, (Y')$ is formed by two components $Y\boxdot_\pm (Y''-Y')$ with $\pm \in \{-,+\}$. These are actually boxes. Together with the integer $n+M+1+n''$, each box forms a pair which might be a piece. If it is the case, the triplet is called \emph{admissible}. 
To prove the main theorem, one needs also to consider such a product when the piece $(Y,n)$ is not a puzzle piece. To this end, we ask the piece $(Y,n)$ to be endowed with a graph transform $T_Y$ from the space of nearly horizontal curves in $Y_\se$ with both endpoints in $\partial^s Y_\se$ into itself, such that in particular, the union of the curves in the image of $T_Y$  contains $f^n(Y)$. Then the triplet of pieces with the graph transform $T_Y$ is \emph{pre-admissible} for a parabolic product $\boxdot_\pm$ if the each of the curves of the range of $T_Y$ are folded into $Y'$  but not too close to $\partial^s Y'$, the same condition on orders holds true and the box obtained is not included in an arc of $W^s(A)$. Furthermore, a graph transform associated to this new piece is defined.\medskip

Starting with the simple puzzle pieces (indexed by $\sA_0$) and using the above operations we can construct new pieces. Such pieces can be (canonically) associated to words with letters in the alphabet $\sA_0$ and $\boxdot_\pm$, with  the concatenation standing for the $\star $ product. We formulate on these symbols some combinatorial rules via an alphabet $\sA\supset \sA_0$. We show that each pre-admissible product respecting these rules is automatically admissible. For instance, one of these rules is that  a parabolic product between $(Y,n)$,$(Y',n')$ and $(Y'',n'')$ must satisfies that $(Y',n')$ and $(Y''n'')$ are $\star$-product of the respectively $n$ and $n+1$ term of a sequence of puzzle pieces $(Y_i,n_i)_i$ satisfying Yoccoz' strongly regular condition $(\diamondsuit )$ stated on p. \pageref{SRcondition}.  \medskip

For every map $f$ close to $h_{a\, 0}$, we can consider the set of all composed operations which are admissible. These form a symbolic set $\sR$ of finite words in the alphabet $\sA$. Each of the words $\sg\in \sR$ is associated to a piece $(Y_\sg, n_\sg)$ and a graph transform $T_\sg$. This is the puzzle structure associated to $f$. 

We can also regard the limit of sequences $(\sg_i)_i$ when $n_{\sg_i}\to \infty$. There are two possible limits.  A sequence $\sc=(\sg_i)$ is  in
$\avv{ \sR}$ if for every $j\ge i$, the word $\sg_j$ starts with $\sg_i$ and is equal to the concatenation of $\sg_j$ with a word in $\sR$. Then the boxes $(Y_{\sg_i})_i$ are nested and $W^s_\sc= \bigcap_i Y_{\sg_i}$ is a Pesin stable manifold. It is a nearly vertical curve with endpoints in $I_\se\times \{-\theta, \theta\}$.  
A sequence  $\st=(\sg_i)_i$ is  in $\arr{ \sR}$ if for every $j\ge i$, the words $\sg_j$ ends with $\sg_i$ and is equal to the concatenation of $\sg_j$ with a  word in $\sR$.
 Then the images of the graph transforms $T_{\sg_i}$ associated to $(Y_{\sg_i}, n_{\sg_i})$ are nested and $\hat W^\st_u= \bigcap_i Im T_{\sg_i}$ is a nearly horizontal curve, whose  endpoints are in $\partial^s Y_\se$ and which contains a Pesin unstable manifold $W^\st_u$.   \medskip

The map $f$ is \emph{strongly regular} if for each $\st\in \arr \sR$, there exist $\sc\in \avv \sR$ such that:
\begin{itemize}
\item[$(\star)$]  the curve $\hat W^\st_u$ is folded by $f^{M+1}| Y_\boxdot$ to a curve tangent to $W^s_\sc$.
\item[$(\diamondsuit )$] the sequence $\sc=\avv \lim \sg_i$ is formed mostly by symbols in $\sA_0$ and the curve $W^s_\sc$ is at most exponentially close to $\partial^s Y_{\sg_i}$ in function of $i$. 
\end{itemize}\medskip

Similarly to the one-dimensional case, each strongly regular map leaves invariant an SRB measure. To handle the parameter selection, we define a cominatorial and arithmetical metric on $\sA^{\Z^-}$. This metric satisfies that the map $\st\in \sR\subset \sA^{\Z^-}\mapsto \hat W^\st_u$ is $\theta$-Lipschitz for the $C^1$-topology of nearly horizontal curves. Also for every $k$ we define a combinatorial map $\pi_k: \sA^{\Z^-}\mapsto \sA^{\Z^-}$ whose image is formed by sequence of $\sA$-symbols which either equal  to $\sA_0$  or encode the rules of parabolic products with pieces of low order in function of $k$. Furthermore, under an induction hypothesis called $k$-great regularity, the set $\arr \sR_k':= \pi_k(\arr \sR)$ is included in $\arr \sR$, and $\arr \sR_k'$ is $\theta^k$-dense in  $ \arr \sR$. The map is  $k+1$-\emph{greatly regular} if each curve in  $(\hat W^\st_u)_{\st \in \arr \sR_k'}$ is folded by $f^{M+1}|Y_\boxdot$ at a piece $Y_\sg$ but not too close to the boundary of $\partial^s Y_\sg$, for a word $\sg \in \sR$ containing mostly symbols with order $\le M/2$. A map which is greatly regular for every $k$ is strongly regular. The abundance of strongly regular maps is shown by bounding the cardinality of parameter components on which $\arr \sR'_k$ is constant, the cardinality of $\arr\sR'_k$, and the length of the interval for which a curve $\hat W^\st_u$ given by a $\st\in \arr \sR_k'$ satisfies the latter folding condition\footnote{Up to here, this subsection  has been changed after the death of J.-C. Yoccoz.}. 
  These combinatorial definitions enable one to follow carefully how the whole structure survives by parameter deformation.

\bibliographystyle{alpha}
\bibliography{references}

\def\polhk#1{\setbox0=\hbox{#1}{\ooalign{\hidewidth
  \lower1.5ex\hbox{`}\hidewidth\crcr\unhbox0}}}
\begin{thebibliography}{RHO80}

\bibitem[AM03]{AM03}
A.~Avila and C.~G. Moreira.
\newblock Statistical properties of unimodal maps: smooth families with
  negative {S}chwarzian derivative.
\newblock {\em Ast\'erisque}, (286):xviii, 81--118, 2003.
\newblock Geometric methods in dynamics. I.

\bibitem[BC85]{BC1}
M.~Benedicks and L.~Carleson.
\newblock On iterations of $1-ax^2$.
\newblock {\em Ann. Math.}, 122:1--25, 1985.

\bibitem[BC91]{BC2}
M.~Benedicks and L.~Carleson.
\newblock The dynamics of the {H}\'enon map.
\newblock {\em Ann. Math.}, 133:73--169, 1991.

\bibitem[BC14]{BeCa13}
P.~Berger and P.~D. Carrasco.
\newblock Non-uniformly hyperbolic diffeomorphisms derived from the standard
  map.
\newblock 329(1):239--262, 2014.

\bibitem[Ber]{B11}
P.~Berger.
\newblock Abundance of non-uniformly hyperbolic {H}\'enon like endomorphisms.
\newblock {\em This volume}.

\bibitem[BP06]{BP06}
L.~Barreira and Y.~Pesin.
\newblock Smooth ergodic theory and nonuniformly hyperbolic dynamics.
\newblock In {\em Handbook of dynamical systems. {V}ol. 1{B}}, pages 57--263.
  Elsevier B. V., Amsterdam, 2006.
\newblock With an appendix by Omri Sarig.

\bibitem[BV01]{BV01}
M.~Benedicks and M.~Viana.
\newblock Solution of the basin problem for {H}\'enon-like attractors.
\newblock {\em Invent. Math.}, 143(2):375--434, 2001.

\bibitem[BY93]{BY}
M.~Benedicks and L.-S. Young.
\newblock Sina\u\i-{B}owen-{R}uelle measures for certain {H}\'enon maps.
\newblock {\em Invent. Math.}, 112(3):541--576, 1993.

\bibitem[Car91]{Ca91}
L.~Carleson.
\newblock The dynamics of non-uniformly hyperbolic systems in two variables.
\newblock In {\em Proceedings of the {I}nternational {C}ongress of
  {M}athematicians, {V}ol.\ {I}, {II} ({K}yoto, 1990)}, pages 1241--1247. Math.
  Soc. Japan, Tokyo, 1991.

\bibitem[CE80]{CE80}
P.~Collet and J.-P. Eckmann.
\newblock On the abundance of aperiodic behaviour for maps on the interval.
\newblock {\em Comm. Math. Phys.}, 73(2):115--160, 1980.

\bibitem[CE83]{CE83}
P.~Collet and J.-P. Eckmann.
\newblock Positive {L}iapunov exponents and absolute continuity for maps of the
  interval.
\newblock {\em Ergodic Theory Dynam. Systems}, 3(1):13--46, 1983.

\bibitem[Dua94]{Du94}
P.~Duarte.
\newblock Plenty of elliptic islands for the standard family of area preserving
  maps.
\newblock {\em Ann. Inst. H. Poincar\'e Anal. Non Lin\'eaire}, 11(4):359--409,
  1994.

\bibitem[G{\'S}97]{GS97}
J.~Graczyk and G.~{\'S}wiatek.
\newblock Generic hyperbolicity in the logistic family.
\newblock {\em Ann. of Math. (2)}, 146(1):1--52, 1997.

\bibitem[H{\'e}n76]{Henon}
M.~H{\'e}non.
\newblock A two-dimensional mapping with a strange attractor.
\newblock {\em Comm. Math. Phys.}, 50(1):69--77, 1976.

\bibitem[Hub93]{Hu91}
J.~H. Hubbard.
\newblock Local connectivity of {J}ulia sets and bifurcation loci: three
  theorems of {J}.-{C}. {Y}occoz.
\newblock In {\em Topological methods in modern mathematics ({S}tony {B}rook,
  {NY}, 1991)}, pages 467--511. Publish or Perish, Houston, TX, 1993.

\bibitem[Jak81]{Ja81}
M.~V. Jakobson.
\newblock Absolutely continuous invariant measures for one-parameter families
  of one-dimensional maps.
\newblock {\em Comm. Math. Phys.}, 81(1):39--88, 1981.

\bibitem[Lor63]{Lorenz}
E.~N. Lorenz.
\newblock Deterministic nonperiodic flow.
\newblock {\em J. Atmos. Sci.}, 20:130--141, 1963.

\bibitem[Luz00]{Lu00}
S.~Luzzatto.
\newblock Bounded recurrence of critical points and {J}akobson's theorem.
\newblock In {\em The {M}andelbrot set, theme and variations}, volume 274 of
  {\em London Math. Soc. Lecture Note Ser.}, pages 173--210. Cambridge Univ.
  Press, Cambridge, 2000.

\bibitem[LY85]{LY}
F.~Ledrappier and L.-S. Young.
\newblock The metric entropy of diffeomorphisms. {I}. {C}haracterization of
  measures satisfying {P}esin's entropy formula.
\newblock {\em Ann. of Math. (2)}, 122(3):509--539, 1985.

\bibitem[Lyu97]{Ly97}
M.~Lyubich.
\newblock Dynamics of quadratic polynomials. {I}, {II}.
\newblock {\em Acta Math.}, 178(2):185--247, 247--297, 1997.

\bibitem[Lyu02]{L02}
M.~Lyubich.
\newblock Almost every real quadratic map is either regular or stochastic.
\newblock {\em Ann. of Math. (2)}, 156(1):1--78, 2002.

\bibitem[Mil00]{Mi92}
J.~Milnor.
\newblock Local connectivity of {J}ulia sets: expository lectures.
\newblock In {\em The {M}andelbrot set, theme and variations}, volume 274 of
  {\em London Math. Soc. Lecture Note Ser.}, pages 67--116. Cambridge Univ.
  Press, Cambridge, 2000.

\bibitem[Mil06]{Mi06}
J.~Milnor.
\newblock {\em Dynamics in one complex variable}, volume 160 of {\em Annals of
  Mathematics Studies}.
\newblock Princeton University Press, Princeton, NJ, third edition, 2006.

\bibitem[MV93]{MV93}
L.~Mora and M.~Viana.
\newblock Abundance of strange attractors.
\newblock {\em Acta Math.}, 171(1):1--71, 1993.

\bibitem[MY10]{MY10}
C.~G. Moreira and J.-C. Yoccoz.
\newblock Tangences homoclines stables pour des ensembles hyperboliques de
  grande dimension fractale.
\newblock {\em Ann. Sci. \'Ec. Norm. Sup\'er. (4)}, 43(1):1--68, 2010.

\bibitem[New74]{Newhouse}
S.~E. Newhouse.
\newblock Diffeomorphisms with infinitely many sinks.
\newblock {\em Topology}, 12:9--18, 1974.

\bibitem[PT93]{PT93}
J.~Palis and F.~Takens.
\newblock {\em Hyperbolicity and sensitive chaotic dynamics at homoclinic
  bifurcations}, volume~35 of {\em Cambridge Studies in Advanced Mathematics}.
\newblock Cambridge University Press, Cambridge, 1993.
\newblock Fractal dimensions and infinitely many attractors.

\bibitem[PY09]{PY09}
J.~Palis and J.-C. Yoccoz.
\newblock Non-uniformly hyperbolic horseshoes arising from bifurcations of
  {P}oincar\'e heteroclinic cycles.
\newblock {\em Publ. Math. Inst. Hautes \'Etudes Sci.}, (110):1--217, 2009.

\bibitem[Ree86]{Re86}
M.~Rees.
\newblock Positive measure sets of ergodic rational maps.
\newblock {\em Ann. Sci. \'Ecole Norm. Sup. (4)}, 19(3):383--407, 1986.

\bibitem[RHO80]{RHO80}
David~A. Russell, James~D. Hanson, and Edward Ott.
\newblock Dimension of strange attractors.
\newblock {\em Phys. Rev. Lett.}, 45(14):1175--1178, 1980.

\bibitem[Ryc88]{R88}
M.~R. Rychlik.
\newblock Another proof of {J}akobson's theorem and related results.
\newblock {\em Ergodic Theory Dynam. Systems}, 8(1):93--109, 1988.

\bibitem[Sar13]{Sa13}
Omri~M. Sarig.
\newblock Symbolic dynamics for surface diffeomorphisms with positive entropy.
\newblock {\em J. Amer. Math. Soc.}, 26(2):341--426, 2013.

\bibitem[Shu71]{Sh71}
M.~Shub.
\newblock {\em Topologically transitive diffeomorphisms of T4}.
\newblock Lecture Notes in Mathematics, Vol. 206. Springer-Verlag, Berlin-New
  York, 1971.

\bibitem[Sin94]{Si94}
Ya.~G. Sina{\u\i}.
\newblock {\em Topics in ergodic theory}, volume~44 of {\em Princeton
  Mathematical Series}.
\newblock Princeton University Press, Princeton, NJ, 1994.

\bibitem[Sma67]{Sm}
S.~Smale.
\newblock Differentiable dynamical systems.
\newblock {\em Bull. Amer. Math. Soc.}, 73:747--817, 1967.

\bibitem[Tak11]{Ta11}
H.~Takahasi.
\newblock Abundance of non-uniform hyperbolicity in bifurcations of surface
  endomorphisms.
\newblock {\em Tokyo J. Math.}, 34(1):53--113, 2011.

\bibitem[Tsu93a]{Ts932}
M.~Tsujii.
\newblock Positive {L}yapunov exponents in families of one-dimensional
  dynamical systems.
\newblock {\em Invent. Math.}, 111(1):113--137, 1993.

\bibitem[Tsu93b]{Ts93}
M.~Tsujii.
\newblock A proof of {B}enedicks-{C}arleson-{J}acobson theorem.
\newblock {\em Tokyo J. Math.}, 16(2):295--310, 1993.

\bibitem[Ure95]{Ur95}
R.~Ures.
\newblock On the approximation of {H}\'enon-like attractors by homoclinic
  tangencies.
\newblock {\em Ergodic Theory Dynam. Systems}, 15(6):1223--1229, 1995.

\bibitem[VL03]{VL03}
M.~Viana and S.~Lutstsatto.
\newblock Exclusions of parameter values in {H}\'{e}non-type systems.
\newblock {\em Uspekhi Mat. Nauk}, 58(6(354)):3--44, 2003.

\bibitem[WY01]{WY01}
Q.~Wang and L.-S. Young.
\newblock Strange attractors with one direction of instability.
\newblock {\em Comm. Math. Phys.}, 218(1):1--97, 2001.

\bibitem[WY08]{YW08}
Q.~Wang and L.-S. Young.
\newblock Toward a theory of rank one attractors.
\newblock {\em Ann. of Math. (2)}, 167(2):349--480, 2008.

\bibitem[Yoc]{Y95}
J.-C. Yoccoz.
\newblock A proof of {J}akobson's theorem.
\newblock {\em This volume}.

\bibitem[You98]{Y98}
L.-S. Young.
\newblock Statistical properties of dynamical systems with some hyperbolicity.
\newblock {\em Ann. of Math. (2)}, 147(3):585--650, 1998.

\end{thebibliography}
\end{otherlanguage}
\end{document}